\documentclass[12pt,leqno]{amsart}
\usepackage{amsmath,stmaryrd}
\usepackage{amssymb}
\usepackage{amsthm}
\usepackage{epsf}
\usepackage{graphicx}
\usepackage{esint} 
\usepackage{dsfont}
\usepackage[pagebackref]{hyperref} 
\hypersetup{pdfpagemode=FullScreen,  colorlinks=true} 
\usepackage{enumerate} 

\numberwithin{equation}{section}

\newcommand{\ms}{\medskip}

\newcommand{\R}{\mathbb{R}}

\newcommand{\B}{\mathcal{B}}
\newcommand{\T}{\mathcal{T}}

\newcommand{\bN}{\mathbb{N}}

\newcommand{\dist}{\,\mathrm{dist}}

\newcommand{\wt}{\widetilde}
\newcommand{\wh}{\widehat}

\newcommand{\1}{{\mathds 1}}

\DeclareMathOperator*{\osc}{osc}

\DeclareMathOperator{\diver}{div}

\theoremstyle{plain}
\newtheorem{theorem}[equation]{Theorem}
\newtheorem{lemma}[equation]{Lemma}
\newtheorem{corollary}[equation]{Corollary}
\newtheorem{proposition}[equation]{Proposition}

\newtheorem{definition}[equation]{Definition}

\theoremstyle{definition}

\theoremstyle{remark}
\newtheorem{example}[equation]{Example}
\newtheorem{remark}[equation]{Remark}

\newcommand{\re}{\mathbb{R}}

\newcommand{\A}{\mathcal{A}}
\newcommand{\C}{\mathcal{C}}

\newcommand{\bp}{\noindent {\it Proof}.\,\,}
\newcommand{\ep}{\hfill$\Box$ \vskip 0.08in}

\def\div{\mathop{\operatorname{div}}}

\begin{document}

\title[]{The Green function with pole at infinity applied to the study of the elliptic measure.}

\author[Feneuil]{Joseph Feneuil}
\address{Joseph Feneuil. Universit\'e Paris-Saclay, Laboratoire de math\'ematiques d’Orsay, 91405, Orsay, France}
\email{joseph.feneuil@universite-paris-saclay.fr}

\begin{abstract}
In $\R^{d+1}_+$ or in $\R^n\setminus \R^d$ ($d<n-1$), we study the Green function with pole at infinity introduced by David, Engelstein, and Mayboroda in \cite{DEM}. In two cases, we deduce the equivalence between the elliptic measure and the Lebesgue measure on $\R^d$; and we further prove the $A_\infty$-absolute continuity of the elliptic measure for operators that can be related to the two previous cases via Carleson measures, extending the range of operators for which the $A_\infty$-absolute continuity of the elliptic measure is known.
\end{abstract}

\maketitle

\ms\noindent{\bf Key words:}
Green function with pole at infinity, elliptic measure, $A_\infty$-absolute continuity.

\ms\noindent

\tableofcontents

\section{Introduction}

\subsection{History and motivation}

Over the past decades, a considerable amount of articles was produced to study the relationship between the geometry of the boundary of a domain $\Omega$ and the $L^p$-solvability of the Dirichlet problem $-\Delta u = 0$ in $\Omega$. The $L^p$-solvability of the Dirichlet problem for large $p$ is equivalent to the absolute continuity of the harmonic measure, and we shall focus our presentation on the later. The theory pioneered with the Riesz brothers in 1916 (see \cite{RR}), who established the absolute continuity of the harmonic measure for simply connected domains in the complex plane with a rectifiable boundary. The quantitative and local analogues are stated in \cite{Lv} and \cite{BJ} respectively. The development of the theory in $\R^n$, $n\geq 2$, started with \cite{Da} and treated Lipschitz domains. Many works were then devoted to find the optimal conditions on $\Omega$ and $\partial \Omega$ to guarantee the absolute continuity of the harmonic measure. It was finally understood that a quantitative version of absolute continuity of the harmonic measure holds if and only if the boundary $\partial \Omega$ is uniformly rectifiable and the domain $\Omega$ has enough access to its boundary. A non-exhaustive list of articles that lead to this conclusion includes \cite{DJ}, \cite{Se}, \cite{HM1}, \cite{HMU}, \cite{AHMNT}, \cite{AHM3TV}, and the minimal access condition to the boundary was recently obtained in \cite{AHMMT}.

One of the strategies to study the absolute continuity of the harmonic measure, and by extension the $L^p$-solvability of the Dirichlet problem, is to make a change of variable in order to obtain an equivalent problem for simpler sets but for more complicated elliptic operators. So instead of studying $-\Delta u = 0$ on a general domain $\Omega$,  many works focused their interest on the study of elliptic operator in the form $L = - \div \A \nabla $ on, for instance, $\Omega = \R^{n-1}_+ := \{(x,t) \in \R^{n-1} \times (0,+\infty)\}$. Here, $\A$ is a matrix satisfying the elliptic and boundedness conditions
 \begin{equation} \label{defelliptic}
 \A(x,t) \xi \cdot \xi  \geq C_L |\xi|^2 \qquad \text{ for } (x,t) \in \Omega \text{ and } \xi \in \R^n,
 \end{equation}
and
 \begin{equation} \label{defbounded}
  |\A(x,t) \xi \cdot \zeta|  \leq C_L |\xi| |\zeta |\qquad \text{ for } (x,t) \in \Omega \text{ and } \xi,\zeta \in \R^n
 \end{equation}
 for some constant $C_L>0$.  As showed in \cite{CFK} and \cite{MM}, the conditions \eqref{defelliptic} and \eqref{defbounded} are not sufficient to ensure that the elliptic measure associated to $L$ is absolutely continuous with respect to the Lebesgue measure on $\R^{n-1}$, and thus some extra assumptions are needed on $\A$ to obtain our absolute continuity.
 Two situations that give positive results are heavily studied: the first situation focuses on $t$-independent matrices $\A$ and are studied in \cite{JK} (use a Rellich identity), \cite{AAH} (perturbations), \cite{HKMP} ($\A$ is non-symmetric), \cite{HLM} (Dirichlet problem in weighted $L^p$), or \cite{HLMP} (the antisymmetric part of $\A$ can be unbounded); while in the second situation, the coefficients of $\A$ satisfy some conditions described with the help of Carleson measures and Carleson measure perturbations, and are considered in for instance \cite{FKP}, \cite{KePiDrift}, \cite{DPP}, \cite{DPR}, \cite{HM12}, \cite{DP},\cite{HMMTZ}, and \cite{HMMTZ}. 
 
When the domain is the complement of a thin set, for instance $\Omega = \R^n \setminus \R^d := \{(x,t) \in \R^d \times \R^{n-d}, \, t \neq 0\}$ with $d < n-1$, studying the solutions to $-\Delta u = 0$ in $\Omega$ does not make sense. 
Indeed, the solutions to $-\Delta u = 0$ in $\Omega$ are the same as the solutions to $-\Delta u = 0$ in $\R^n$, which means that the boundary $\R^d$ is not ``seen'' by the Laplacian or, in term of harmonic measure, it means that the Brownian motion has zero probability to hit the boundary $\R^d$. 
G. David, S. Mayboroda, and the author developed in \cite{DFMprelim}, \cite{DFMprelim2} an elliptic theory for domains with thin boundaries by using appropriate degenerate operators. If $\Omega = \R^n \setminus \R^d$ is considered, we assume that the elliptic operator $L=-\div A \nabla$ satisfies
 \begin{equation} \label{defelliptic2}
 A(x,t) \xi \cdot \xi  \geq C_L |t|^{d+1-n} |\xi|^2 \qquad \text{ for } (x,t) \in \Omega \text{ and } \xi \in \R^n,
 \end{equation}
and
 \begin{equation} \label{defbounded2}
  |A(x,t) \xi \cdot \zeta|  \leq C_L |t|^{d+1-n} |\xi| |\zeta |\qquad \text{ for } (x,t) \in \Omega \text{ and } \xi,\zeta \in \R^n,
 \end{equation}
for some constant $C_L >0$. The operator $L$ can thus be written $-\div |t|^{d+1-n} \A \nabla$ and $\A$ satisfies \eqref{defelliptic}--\eqref{defbounded}. Under those conditions, the elliptic measure with pole in $X\in \Omega$ associated to $L$, denoted by $\omega_L^X$, is the probability measure on $\R^d$ so that the function $u_f$ on $\Omega$ constructed for any $f\in C^\infty_0(\R^d)$ as
\begin{equation} \label{defhm}
u_f(X) = \int_{\R^d} f(y) d\omega_L^X(y)
\end{equation}
is a weak solution to $Lu=0$, continuous on $\overline\Omega$, and whose trace on $\R^d$ is $f$.
The articles \cite{DFMAinfty}, \cite{MZ}, \cite{FMZ}, \cite{DM}, \cite{Fen}, \cite{MP} tackled the absolute continuity of the elliptic measure (or $L^p$-solvability of Dirichlet problem) in the case where the boundary of $\Omega$ is a low dimensional set.

We finish the subsection with the following observation done in \cite{DFMKenig}. 
Let $L= - \div \A \nabla$ be an elliptic operator defined on $\R^{d+1}_+$ and which satisfies  \eqref{defelliptic}--\eqref{defbounded}. We define $\A_{1}$ as the $d\times d$ submatrix of $\A$ in the top-left corner, and $\A_2$, $\A_3$, $a_4$ so that $\A$ is written by block as
\begin{equation} \label{decompA}
\A = \begin{pmatrix} 
\A_1 & \A_2 \\
\A_3 & a_4
\end{pmatrix}
\end{equation}
For $n>d+1$, we construct the elliptic operator $\wt L= - \div |t|^{d+1-n} \wt \A \nabla$ defined on $\R^n\setminus \R^d$ as
\begin{equation} \label{defwtA}
\wt \A(x,t) := \begin{pmatrix} 
\A_1(x,|t|)  & \A_2(x,|t|) \frac{t}{|t|} \\
\frac{t^T}{|t|} \A_3(x,|t|) & a_4(x,|t|) I_{n-d}
\end{pmatrix},
\end{equation}
where $t$ is seen here as a horizontal vector in $\R^{n-d}$, which means that $\A_2t$ and $t^T\A_3$ are matrices of dimensions $d\times (n-d)$ and $(n-d) \times d$ respectively, and $I_{n-d}$ is the identity matrix of order $n-d$.  Then the elliptic measures on $\R^d$ associated to $L$ and $\wt L$ - we call them $\omega^{(x,r)}$ and $\wt \omega^{(x,t)}$ - satisfy that
\begin{equation}
\wt \omega^{(x,t)} = \omega^{(x,|t|)} \qquad \text{ for $(x,t) \in \R^n\setminus \R^d$.}
\end{equation}
More generally, any solution $u$ to $Lu = 0$ in $\R^{d+1}_+$ yields a solution $\wt u(x,t) := u(x,|t|)$ to $\wt Lu = 0$ in $\R^n \setminus \R^d$. As a consequence, the construction from \cite{CFK,MM} can be adapted to provide, for any $1\leq d <n$, examples of operators whose elliptic measures are not absolutely continuous with respect to the Lebesgue measure on $\R^d$. 
It also means that if an operator $\wt L = - \div |t|^{d+1-n} \wt \A \nabla$ can be written as \eqref{defwtA}, and if the elliptic measure of the original operator $L=-\div \A \nabla$ is absolutely continuous with respect to the Lebesgue measure on $\R^d$, then the elliptic measure associated to $\wt L$ is also absolutely continuous with respect to the Lebesgue measure on $\R^d$. 
The above construction provides, for any dimension and codimension of the boundary, a wide range of elliptic operators that satisfy the absolute continuity of the elliptic measure. However, the (relevant) solutions of those operators are radial, i.e. depends only on the distance to the boundary $\R^d$ and its projection on $\R^d$. 

The goal of this article is to go beyond the matrices that can be written as \eqref{defwtA}. Of course, as we shall discuss in the next subsection, we already know that some cases where the first $d$ lines do not matter for the $A_\infty$-absolute continuity of the elliptic measure (see \cite{DFMAinfty}, \cite{FMZ}), and we also know that the $A_\infty$ is stable under Carleson perturbations (see \cite{MP}). However, we do not know, for instance, whether it is possible that the bottom right corner of $\A$ is not a Carleson perturbation of a sub-matrix in the form $b(x,|t|)I_{n-d}$.

In this article, we will construct the Green function with pole at infinity, in the spirit of \cite{DEM}. The Green function (and the Green function with pole at infinity) has a deep connection with the harmonic measure (see Lemma \ref{lemG=tw} below); some recents works (\cite{DM2}, \cite{DFM-GFE}, and \cite{DLM}) even started to link the geometry of $\partial \Omega$ directly to bounds on the Green function (instead of estimates on the harmonic measure). We shall thus study the Green function with pole at infinity in few easy cases, and deduce that the harmonic measure and the Lebesgue measure are comparable (hence $A_\infty$ absolute continuous with respect to each other). Then we shall use the Green function with pole at infinity as a substitute of $|t|$ in a now classical argument that establish the stability of the $A_\infty$-absolute continuity of the harmonic measure under some transformations on the elliptic operator. This will enlarge the class of operators for which the $A_\infty$-absolute continuity of the harmonic measure is known, and especially in the case where $d<n-1$.

\subsection{Presentation of the results}

In the rest of the article $d$ is an integer in $\{0,\dots,n-1\}$. If $d=n-1$, then $\Omega = \R^{n}_+ = \R^{d+1}_+ = \{(x,t) \in \R^d \times (0,+\infty)\}$. 
If $d<n-1$, then $\Omega = \R^n \setminus \R^d = \{(x,t) \in \R^d \times \R^{n-d}, \, t \neq 0 \}$. When we write that $0< |t|< r$, we understand $t\in (0,r)$ if $n-d=1$ and $t \in B(0,r) \subset \R^{n-d}$ otherwise.

 If $L=-\div A \nabla$ satisfies \eqref{defelliptic2}--\eqref{defbounded2}, then the elliptic measure defined in \eqref{defhm} is non-degenerate, is doubling, and satisfies the change of pole property (respectively Lemmas 11.10, 11.12, and 11.16 in \cite{DFMprelim}), and those conditions are the ones needed to prove the following result from \cite{DFMAinfty}.

\begin{theorem}[Theorem 8.9 in \cite{DFMAinfty}] \label{ThAinfty}
Let $L = - \div A \nabla$, where the real matrix-valued function $A$ satisfies the ellipticity and boundedness conditions \eqref{defelliptic}--\eqref{defbounded}. Assume that there exists $M>0$ such that, for any Borel set $H\subset \R^d$, the solution $u_H$ defined by $u_H(X) = \omega_L^X(H)$ satisfies the Carleson measure estimate
\begin{equation} \label{NablauCM}
\sup_{x\in \R^d, \, r>0} \fint_{B_{\R^d}(x,r)} \int_{0< |t| < r} |t\nabla u_H|^2 \frac{dy \, dt}{|t|^{n-d}} \leq M.
\end{equation}
Then the elliptic measure is $A_\infty$ with respect to the Lebesgue measure on $\R^d$, i.e. for every $\epsilon>0$, there exists a $\delta >0$ (that depends only on $\epsilon$, $d$, $n$, $C_L$, and $M$) such that for every ball $B:=B(x,r) \subset \R^d$, every $t$ that verifies $|t|=r$, and any Borel set $E \subset B$, one has
\begin{equation}
\text{ if } \frac{\omega_L^{(x,t)}(E)}{\omega_L^{(x,t)}(B)} < \delta, \text{ then } \frac{|E|}{|B|} < \epsilon.
\end{equation}
\end{theorem}

For a proof when $d=n-1$, see Corollary 3.2 in \cite{KKiPT}. The condition \eqref{NablauCM} is closely related to another characterization for the $A_\infty$-absolute continuity of the elliptic measure on $\R^d$ called $BMO$-solvability, which can be found in \cite{DKP} for the codimension 1 and in \cite{MZ} when $d<n-1$. 

The condition \eqref{NablauCM} means that $(|t||\nabla u_H|)^2 |t|^{d-n} dt\, dx$ is a Carleson measure. To lighten the presentation, we even introduce a notation for inequalities like \eqref{NablauCM}. We say that a quantity $f$ satisfies the {\bf Carleson measure condition} if there exists $C>0$ such that 
\begin{equation} \label{defCM}
\|f\|_{L^\infty} \leq C \quad \text{ and } \quad \sup_{x\in \R^d, \, r>0} \fint_{B_{\R^d}(x,r)} \int_{|t| < r} |f|^2 \frac{dt \, dy}{|t|^{n-d}} \leq C.
\end{equation}
For short, we write $f \in CM_2$ or $f\in CM_2(C)$ when we want to refer to the constant in the right side of the bound \eqref{NablauCM}. So to conclude, in order to apply Theorem \ref{ThAinfty}, we need to assume that there exists $K>0$ such that for any Borel set $H$, the function $u_H \in CM_2(K)$. It will also be useful to write the variant $f\in \wt{CM}_2(C)$ when 
\begin{equation} \label{defwtCM}
\sup_{x\in \R^d, \, r>0} \fint_{B_{\R^d}(x,r)} \int_{|t| <r} \Big(\sup_{|Z-(y,t)| < |t|/4} |f(Z)|^2\Big) \frac{dt \, dy}{|t|^{n-d}} \leq C.
\end{equation}

\medskip

To the best of the author's knowledge, in our setting of high co-dimensional boundaries, the most general condition on the coefficients of the matrix $A$ that ensures the $A_\infty$ absolute continuity of the elliptic measure with respect to the $d$-dimensional Hausdorff measure is given in \cite{FMZ}. 

\begin{theorem}[Theorem 1.9 (1) in \cite{FMZ} for $p=2$] \label{ThFMZ}
Let $L = - \div |t|^{d+1-n} \A \nabla$, where the real matrix-valued function $\A$ satisfies the ellipticity and boundedness conditions \eqref{defelliptic}--\eqref{defbounded}. 
Assume that $\A$ can be decomposed as
\begin{equation} \label{A=B+C}
\A = \begin{pmatrix} \A_1 & \A_2 \\ \B_3 & b.I_{n-d} \end{pmatrix} + \C
\end{equation}
with $I_{n-d}$ being the identity matrix, $\A_1$, $\A_2$, and $\B_3$ being respectively $d\times d$, $d\times (n-d)$, and $(n-d) \times d$ matrix-valued functions, $b$ being a scalar function, who satisfies 
\begin{itemize}
\item $K^{-1} \leq b \leq K$,
\item $|t||\nabla b| + |t||\nabla_x \B_3| + |t|^{n-d} \div_t(|t|^{d+1-n} \B_3) + |\C| \in CM_2(K)$
\end{itemize}
for a constant $K>0$. Then the hypothesis \eqref{NablauCM} of Theorem \ref{ThAinfty} is satisfied (with a constant $M$ that depends only on $d$, $n$, $C_L$, and $K$) and therefore the elliptic measure $\omega_L^X$ is $A_\infty$-absolutely continuous with respect to the Lebesgue measure.
\end{theorem} 

\noindent {\em Remarks}
\begin{enumerate}[(i)]
\item In codimension 1, that is when $d=n-1$, Theorem \ref{ThFMZ} requires that the last line $\mathfrak a_{d+1}$ of the matrix $\A$ can be decomposed as $a_{d+1} = \mathfrak b_{d+1} + \mathfrak c_{d+1}$ with $|t||\nabla \mathfrak b_{d+1}| + |\mathfrak c_{d+1}| \in CM_2$. This condition is thus weaker than the one found in \cite{KePiDrift}, where one assume that $|t||\nabla \A| \in CM$, and the conditions are the same if we add to \cite{KePiDrift} the perturbation theory from \cite{HM12}. However, to the best of the author's knowledge, the first time where no conditions on the first $d$ lines were assumed is in \cite{DFMAinfty} and \cite{FMZ}.
\item Observe that if $\A$ is a $(d+1)\times (d+1)$-matrix valued function on $\R^{d+1}_+$ that satisfies the assumptions of the above theorem, then the $n\times n$-matrix valued function $\wt \A$ defined from $\A$ on $\R^n \setminus \R^d$ as in \eqref{defwtA} satisfies the assumptions of Theorem \ref{ThFMZ} too. 
\item With the same argument as the one used in Corollary 2.3 in \cite{DPP}, one can show that if $\B_3$ and $b$ satisfy
\begin{equation} \label{oscB3b}
(x,t) \mapsto \osc_{B((x,t),|t|/4)} \B_3 + \osc_{B((x,t),|t|/4)} b \in CM_2(K),
\end{equation}
where $\osc_B f = \sup_B f - \inf_B f$, then we can find $\wh{\B_3}$ and $\wh{b}$ such that 
\[(x,t) \mapsto \sup_{B((x,t),|t|/4)} |\B_3 - \wh {\B_3}| +  \osc_{B((x,t),|t|/4)} |b - \wh b| \in CM_2(K')\]
and
\[ \nabla\wh{\B_3}  + \nabla\wh b \in CM_2(K').\]
So assuming the apparently weaker condition \eqref{oscB3b} is enough to satisfy the assumptions of Theorem \ref{ThFMZ} and therefore obtain the $A_\infty$-absolute continuity of the elliptic measure.
\end{enumerate}

When $d<n-1$, the operator $L = -\div A \nabla$ will necessarily depends on $|t|$ as long as it satisfies \eqref{defelliptic2}--\eqref{defbounded2}. However, once the weight $|t|^{d+1-n}$ is removed, we can see that Theorem \ref{ThFMZ} does not even consider the simple case where $\A = |t|^{n-d-1} A$ is an arbitrary constant matrix.

\medskip

Let $\T_3$ be a $(n-d)\times d$ matrix-valued function and $\T_4$ be a $(n-d)\times(n-d)$ matrix-valued function. We say that $(\T_3,\T_4)$ satisfies \eqref{H1} if
\begin{equation} \label{H1} \tag{H1}
\text{$\T_3$ and $\T_4$ are $x$-independent.}
\end{equation} 
And we say that $(\T_3,\T_4)$ satisfies \eqref{H2} if
\begin{equation} \label{H2} \tag{H2}
\left\{\begin{array}{l}
\text{$(\T_3)^T \nabla |t|$ is $x$-independent} \\
\text{and there exists $h: \, (0,+\infty) \mapsto \R$ such that $(\T_4)^T \nabla |t| = h(|t|) \nabla |t|$.}
\end{array}\right.
\end{equation}
In addition, we say that $\T_4$ satisfies \eqref{H1}/\eqref{H2} if $(0,\T_4)$ satisfies \eqref{H1}/\eqref{H2}. Note that when $d=n-1$, $\T_4$ is a scalar function, and \eqref{H1} and \eqref{H2} are the same hypothesis.

The condition \eqref{H2} for $\T_4$ is neither weaker nor sronger than \eqref{H1}. Roughly, if $\T_4$ safisties \eqref{H2}, then $\nabla |t|$ is an eigenvalue of $\T_4$ and $\T_4$ may depend of $x$. 

\begin{example} \label{exT4}
If we set $v(t)$ to be an horizontal vector orthogonal to $t$ and independent of $x$, for instance $v(t) = (-t_2,t_1,0 \dots,0)$, and $a(x)$ is an vertical vector in $\R^{n-d}$ independent of $t$, for instance $a(x) = (\cos(x),0,\dots,0)^T$, then 
\[\T_4:= I_{n-d} + \frac1{2|t|} a(x)v(t)\]
satisfies \eqref{H2} but not \eqref{H1}. On the opposite, a matrix $\T_4$ which is constant will satisfy \eqref{H1} but not \eqref{H2} except if $\T_4$ is actually a scalar multiplication of the identity matrix. Also, observe that $\T_4$ can go beyond $b.I_{n-d}$ and still stabilizes $\nabla |t|$. Remember that $t$ is seen as an horizontal vector and hence
\[\T_4:= I_{n-d} + \frac12 \frac{t^Tt}{|t|^2}\]
is a matrix that verifies both \eqref{H1} and \eqref{H2} but is not the multiplication of the identity matrix by a scalar function. 
\end{example}

\medskip

Our first result states that, if the last $n-d$ lines of $\A$ satisfy either \eqref{H1} or \eqref{H2}, then the elliptic measure and the Lebesgue measure on $\R^d$ are equivalent. Taking matrices as given in Example \ref{exT4} will already allow us to obtain control of the harmonic measure for elliptic operators that goes beyond the ones found in the literature. 

\begin{theorem} \label{Main2} Let $L=-\div |t|^{d+1-n} \A \nabla$ be an elliptic operator satisfying \eqref{defelliptic}--\eqref{defbounded}. Assume that $L$ is such that 
\begin{equation} \label{A=AT}
\A = \begin{pmatrix} \A_1 & \A_2 \\ \T_3 & \T_4 \end{pmatrix},  \text{ where $(\T_3,\T_4)$ satisfies either \eqref{H1} or \eqref{H2}.}
\end{equation}
Then, for any $Y_0= (y_0,t_0)$ and any Borel set $E \subset \Delta_{Y_0} := B_{\R^d}(y_0,|t_0|)$, we have
\begin{equation}
C^{-1} \frac{|E|}{|\Delta_{Y_0}|} \leq \omega^{Y_0}(E) \leq \frac{|E|}{|\Delta_{Y_0}|},
\end{equation}
where $|E|$ denotes the $d$-dimensional Lebesgue measure and $C>0$ depends only on $n$, $d$, and $C_L$.
\end{theorem}

For our second result, we consider elliptic operators whose coefficients are close to a matrix in the form \eqref{A=AT}.
We shall show that for such operators the bound \eqref{ThAinfty} holds by adapting a $S<N$ argument (see \cite{KKPT} and ensuing literature). Our contribution will be the use of Green function as a substitute for $|t|$, a bit like in \cite{AMHT} but we handle the (possible) roughness of the Green function in a much simpler Cacciopoli-type argument.

\begin{theorem} \label{Main1}
Let $L=-\div |t|^{d+1-n} \A \nabla$ be an elliptic operator satisfying \eqref{defelliptic}--\eqref{defbounded}. Decompose $\A$ as 
\begin{equation} \label{A=B+C3}
\A = \begin{pmatrix} \A_1 & \A_2 \\ \B_3 + \C_3 & b \T_4 + \C_4 \end{pmatrix}, 
\end{equation}
where $b$ is a scalar function, $\A_1$ is a $d\times d$ matrix, and the dimensions of $\A_2$, $\B_3$, $\C_3$, $\T_4$, $\C_4$ are such that the matrices complete the $n\times n$-matrix $\A$. Assume that the sub-matrices of $\A$ satisfy that  
\begin{enumerate}[\em (a)]
\item $\T_4$ satisfies either \eqref{H1} or \eqref{H2},
\end{enumerate}
and there exists a constant $K>0$ such that
\begin{enumerate}[\em (a)]
\addtocounter{enumi}{1}
\item $K^{-1} \leq b \leq K$,
\item  $|\C_3| + |\C_4| \in \wt{CM}_2(K)$,
\item $|t||\nabla b| + |t||\div_x (\B_3)^T| + |t|^{n-d} |\div_t (|t|^{d+1-n} \B_3)| \in CM_2(K)$,
\end{enumerate}
Then the hypothesis \eqref{NablauCM} of Theorem \ref{ThAinfty} is true and thus the elliptic measure $\omega^X_L$ is $A_\infty$-absolutely continuous with respect to the Lebesgue measure on $\R^d$.
\end{theorem}

\noindent {\em Remarks:}
\begin{enumerate}[(i)]
\item Theorem \ref{ThFMZ} is a consequence of Theorem \ref{Main1} when  $\T_4$ is the identity matrix. 

\smallskip

\item In the above theorem, when $M = (M_{ij})_{ij}$ is a $(n-d) \times d$-matrix, then the quantity $\div_x M^T$ is a vector in $\R^{n-d}$ whose $k^{th}$ value is $\sum_{j=1}^d \partial_{x_j} M_{kj}$, and similarly the quantity $\div_t M$ is a vector in $\R^d$ whose $k^{th}$ value is $\sum_{i=1}^{n-d} \partial_{t_i} M_{ik}$.

\smallskip

\item When $d=n-1$, $\T_4$ is a scalar function, and \eqref{H2} should read `there exists $h:\, (0,+\infty) \mapsto \R$ such that $\T_4 \nabla t  = h(t) \nabla t$ for all $t\in (0,+\infty)$', but the later just means that $\T_4$ is $x$-independent, and thus \eqref{H1} and \eqref{H2} are the same hypothesis.

\smallskip

\item We actually prove a stronger estimate than \eqref{NablauCM}, we prove a local $S<N$ $L^2$-estimate, which is stated in \eqref{S<N} below. We see {\em a priori} no big obstacles in our methods that will stop us from obtaining $N<S$ estimates under the assumptions of Theorem \ref{Main1}, and hence from studing the solvability of the Dirichlet, Neumann, and regularity problems.
\end{enumerate}

In the next results, we shall assume a stronger condition on $\T_4$, which will allow us to be slightly more flexible on the bottom left corner of $\A$. In the next lemma, $\B_3$ can satisfy either $|t|^{n-d} |\div_t |t|^{d+1-n} \B_3| \in CM_2$ as in Theorem \ref{Main1}, or simply $|t||\div_t \B_3| \in CM_2$.

\begin{theorem} \label{Main3}
Assume that $d < n-2$. Let $L=-\div |t|^{d+1-n} \A \nabla$ be an elliptic operator satisfying \eqref{defelliptic}--\eqref{defbounded}. Decompose $\A$ as 
\begin{equation} \label{A=B+C4}
\A = \begin{pmatrix} \A_1 & \A_2 \\ \B_3 + \wt \B_3 + \C_3 & b \T_4 + \C_4 \end{pmatrix}, 
\end{equation}
and assume that 
\begin{enumerate}[\em (a)]
\item $(\T_4)^T \nabla |t| = \nabla |t|$,
\end{enumerate}
and there exists a constant $K>0$ such that
\begin{enumerate}[\em (a)]
\addtocounter{enumi}{1}
\item $K^{-1} \leq b \leq K$,
\item $|\C_3| + |\C_4| \in \wt{CM}_2(K)$,
\item $|t||\nabla b| + |t||\div_x (\B_3)^T| + |t|^{n-d} |\diver_t (|t|^{d+1-n} \B_3)| \in CM_2(K)$,
\item $|t||\div_x (\wt \B_3)^T| + |t| |\div_t \wt\B_3| \in CM_2(K)$.
\end{enumerate}
Then the hypothesis \eqref{NablauCM} of Theorem \ref{ThAinfty} is true and thus the elliptic measure $\omega^X_L$ is $A_\infty$-absolutely continuous with respect to the Lebesgue measure on $\R^d$.
\end{theorem}

\begin{remark}
The last theorem is a bit unmotivated at the moment. One classical strategy to deal with non-flat boundaries is to make changes of variable. One can thus obtain an equivalent problem where the boundary is better (i.e. flat) but the coefficients of the operators are much worse. See for instance \cite{KePiDrift} in the classical case and \cite{DFMAinfty} in higher codimension. That is why it is key to obtain the largest set of operator that we can understand in the flat case. The term $\B_3$ is the one that we can treat if we adapt the proof of \cite{KePiDrift} in higher codimension, however $\B_3$ is not well adapted to change of variable and we would prefer $\wt \B_3$ much better.

In \cite{DFMAinfty}, the authors had to introduce a new (and more complicated) change of variable in order to deal with the case where the boundary is the graph of a Lipschitz function. Still the construction is limited to graphs of Lipschitz function with small Lipschitz constant. I claim here that we can deal with big Lipschitz constant if we can allow terms in the form $\wt \B_3$ in the bottom left corner of $\A$, as we do in Theorem \ref{Main3}.

The full construction of the change of variable that maps the graph of an arbitrary Lipschitz function $\varphi: \, \R^d \to \R^{n-d}$ to $\R^d$ and that turns the elliptic operator from \cite{DFMAinfty} into one in the form of \eqref{A=B+C4} will not be done here, since it would be too long and technical (and we don't have a new result to prove with it). We will only give a rough idea via an example. If the Lipschitz function is 
\[\varphi: \, x \in \R \mapsto (ax,0) \in \R^2,\]
then the change of variable that send $\R$ on the graph of $\varphi$ constructed in \cite{KePiDrift} would be
\[\rho_1(x,t_1,t_2) = (x,ax+t_1,t_2),\]
while the one in \cite{DFMAinfty} would be
\[\rho_2(x,t_1,t_2) = (x-c_1t_1,ax+c_2t_1,t_2)\]
where $c_1 = a/\sqrt{1+a^2}$ and $c_2 = 1/\sqrt{1+a^2}$.
If we consider the operator $L = -\diver \delta(X)^{-1} \nabla$, where $\delta(X)$ is the distance between $X$ and the graph of $\varphi$, then using the change of variable
\[\rho_3(x,t_1,t_2) = (x,ax+ct_1,t_2), \text{ with } c= \sqrt{1+a^2}\] will turn $L$ into $L_3 = -\diver |t|^{-1} \A^3 \nabla$ where
\[\A^3 = \begin{pmatrix} c^{-4} & -ac & 0 \\ -ac & 1 & 0 \\ 0 & 0 & 1 \end{pmatrix}\]
is in the form \eqref{A=B+C4}.

\end{remark}

The article is divided as follows. Section \ref{SG} introduces the notion of Green function with pole at infinity, and deduce a relation between this Green function and the elliptic measure that holds whenever $L$ satisfies the elliptic and boundedness conditions \eqref{defelliptic}--\eqref{defbounded}. 
Section \ref{S3} is devoted to the study of operators in the form \eqref{A=AT}, and proves Theorem \ref{Main2}. In Section \ref{S4}, we demonstrate Theorem \ref{Main1} and \ref{Main3} by establishing a local $S<N$ estimate that implies \eqref{NablauCM}.

In the rest of the article, $A \lesssim B$ means that $A\leq CB$ for a constant $C$ whose dependence in the parameters will be recalled or obvious from context. In addition, $A\approx B$ means $A\lesssim B$ and $B \lesssim A$.

\section{General results on the Green function  with pole at infinity} \label{SG}

In all this section, we consider an elliptic operator $L=-\div |t|^{d+1-n} \A \nabla$ satisfying \eqref{defelliptic}--\eqref{defbounded}. Even if this article is written when $\Omega = \R^{d+1}_+$ (if $d=n-1)$ or $\Omega = \R^n \setminus \R^d$ (if $d<n-1)$, the definitions and results of this section can be easily generalized to domains (and elliptic operators) that enter the scope of the elliptic theory developed in \cite{DFMprelim} and \cite{DFMprelim2}. In particular, we only need $\Omega$ to satisfy the Harnack chain condition and the corkscrew point condition (see \cite{DFMprelim,DFMprelim2} for the definitions).

We need a bit of functional theory, which is only needed for the precise statement of Definition \ref{defGinfty} and Proposition \ref{propGinfty} below, and can be overlooked. The space
\begin{equation} \label{defW}
W := \{ u \in L^1_{loc}(\Omega) \, \int_\Omega |\nabla u| \frac{dt \, dx}{|t|^{n-d-1}} < +\infty \}
\end{equation}
is equipped with the semi norm $\|u\|_W := \|\nabla u\|_{L^2(\Omega)}$. Observe that $\|.\|_W$ is a norm for $C^\infty_0(\Omega)$ and we write $W_0$ for the completion of $C^\infty_0(\Omega)$ under $\|.\|_W$. We also define 
\begin{equation} \label{defW0O}
W_0(\overline{\Omega}) := \{ u \in W^{1,2}_{loc}(\Omega), \, u\varphi  \in W_0 \text{ for any } \varphi \in C^\infty_0(\R^n) \,  < +\infty \}.
\end{equation}
The proof of the properties of $W$, $W_0$, and $W_0(\overline{\Omega})$ can be found in \cite{DFMprelim} and \cite{DFMprelim2}, but let us give few comments to build the reader's intuition. 
The spaces $W$ and $W_0$ are the ones where we find the solutions to the Dirichlet problem $Lu = 0$ in $\Omega$, $u=f \in H^{1/2}(\R^d)$ by using the Lax-Milgram theorem; here $H^{1/2}(\R^d) = W^{1/2,2}(\R^d) = B^{1/2}_{2,2}(\R^d)$ is the (classical) Besov space of traces. The space $W_0$ is the subspace of $W$ containing the functions with zero trace. The space $W_0(\overline{\Omega})$ is a space bigger than $W_0$, that possess the same local properties as $W_0$, but does not have any control when $|(x,t)| \to \infty$.

We recall that $u\in W^{1,2}_{loc}(\Omega)$ is a weak solution to $Lu= 0$ in $\Omega$ if 
\begin{equation} \label{defsol}
\int_\Omega \A \nabla u \cdot \nabla \varphi \frac{dt\, dx}{|t|^{n-d-1}} = 0 \qquad \text{ for } \varphi \in C^\infty_0(\Omega).
\end{equation}

\begin{definition} \label{defGinfty}
A Green function (associated to $L^*$) with pole at infinity is a {\bf positive} weak solution $G := G_{L^*} \in W_0(\overline{\Omega})$ to $L^*u =-\div |t|^{d+1-n} \A^T \nabla u = 0$ in $\Omega$.
\end{definition}

Be aware that, in the above definition, the function $G$ is a solution to the adjoint operator $L^* = -\div |t|^{d+1-n} \A^T \nabla$. We prefer here to associate $G$ to the adjoint right away, because it is the appropriate tool we ultimately need for our proofs. But since $L$ and $L^*$ satisfy the same properties \eqref{defelliptic}--\eqref{defbounded}, we have the following:

\begin{proposition}[Lemma 6.5 in \cite{DEM}] \label{propGinfty}
A Green function with pole at infinity $G$ enjoys the following properties:
\begin{itemize}
\item $G \in C(\overline{\Omega})$, i.e. $G$ is continuous up to the boundary $\R^d$,
\item $G = 0$ on $\R^d$,
\item $G$ is unique up to a constant. We write $G_{X}$ for the only Green function with pole at infinity which satisfies $G_X(X) = 1$, and the uniqueness gives that
\begin{equation} \label{Gunique}
G_X(Y) G_Y(X) = 1 \quad \text{ for } X,Y \in \Omega.
\end{equation}
\item Let $G^Y(X)$ be the Green function (associated to $L^*$) with pole at $Y$ as defined in Chapter 10 of \cite{DFMprelim}. Take $Y_0 = (y_0,t_0) \in \Omega$, and define for $j\in \bN$ the point $Y_j = (y_0,2^jt_0)$. There exists a subsequence $j_n \to \infty$ such that 
\begin{equation} \label{defG}
\frac{G^{Y_{j_n}}}{G^{Y_{j_n}}(Y_0)} \text{ converges uniformly on compact sets of $\overline \Omega$ to } G_{Y_0}.
\end{equation}
\end{itemize}
\end{proposition}

\bp The first two points are a consequence of the De Giorgi-Nash-Moser estimates on weak solutions that can be found (for instance) in \cite[Chapter 8]{DFMprelim}. The last two points are in Lemma 6.5 from \cite{DEM} or in its proof.
\ep

We assign to any point $Z = (z,s) \in \Omega$, the boundary ball
\begin{equation} \label{defDz}
\Delta_Z := B(z,|s|) \subset \R^d.
\end{equation}
We compare the Green function with pole at infinity and the elliptic measure.

\begin{lemma} \label{lemG=tw}
Let $Y_0 = (y_0,t_0) \in \Omega$. If $X = (x,t) \in \Omega$ verifies $x \in B(y_0,2|t_0|)$ and $0 < |t| < 2|t_0|$, we have that
\begin{equation}
C^{-1} G_{Y_0}(X) \leq \left(\frac{|t|}{|t_0|}\right)^{1-d} \omega^{Y_0}(\Delta_X) \leq C  G_{Y_0}(X),
\end{equation}
where $C>0$ depends only on $n$, $d$, and $C_L$. Here $G_{Y_0} = G_{L^*,Y_0}$ is defined in Proposition \ref{propGinfty} and is the Green function associated to $L^*$ with pole at infinity, and $\omega^{Y_0} = \omega^{Y_0}_L$ is the elliptic measure associated to $L$ with pole at $Y_0$
 defined in \eqref{defhm}.
\end{lemma}

\begin{remark}
We can use the uniqueness of the Green function \eqref{Gunique} to get an estimate of $G_{Y_0}(X)$ using the elliptic measure when $X$ is far from $Y_0$.
\end{remark}

\bp We need to invoke some results from \cite{DFMprelim}. The classical case $d=n-1$ is not included in \cite{DFMprelim}, but is either known from the reader or can be found in the last section of \cite{DFMprelim2}.

Let $Y_1= (y_0,4t_0)$. The change of pole property (Lemma 11.16 in \cite{DFMprelim}) states that, for any Borel set $E \subset \Delta_{Y_1} = 4\Delta_{Y_0}$ and any $Y \in \Omega$ satisfying $|Y-y_0| > 8t_0$, we have that
\begin{equation} \label{chgpole2}
 \omega^{Y_1}(E) \approx \frac{\omega^Y(E)}{\omega^Y(\Delta_{Y_1})}
\end{equation}
with constants that depend only on $n$, $d$, and $C_L$. Together with the doubling property of the elliptic measure (Lemma 11.12 in \cite{DFMprelim}) and the Harnack inequality (Lemma 8.9 in \cite{DFMprelim}), we deduce that, for the same set $E$, point $Y$, and with constants that depend on the same parameters, we have that
\begin{equation} \label{chgpole}
 \omega^{Y_0}(E) \approx \frac{\omega^Y(E)}{\omega^Y(\Delta_{Y_0})}.
\end{equation}

For our second result, we want to compare the Green function and the elliptic measure. Let $g^X(Y)$ be the Green function associated to $L$ with pole in $X$. Then Lemma 10.6 in \cite{DFMprelim} entails that
\begin{equation} \label{compGg}
G^Y(X) = g^X(Y) \qquad \text{ for } X,Y\in \Omega.
\end{equation}
Moreover, Lemma 11.11 in \cite{DFMprelim} gives, for $X=(x,t) \in \Omega$ and $Y \in \Omega \setminus B_{\R^n}(x,2|t|)$, that
\begin{equation} \label{compgw}
|t|^{d-1} g^X(Y) \approx \omega^Y(\Delta_X),
\end{equation}
with constants that depend only on $n$, $d$, and $C_L$. So the combination of \eqref{compGg} and \eqref{compgw} implies, for $X = (x,t) \in \Omega$, that
\begin{equation} \label{compGw}
|t|^{d-1} G^Y(X) \approx \omega^Y(\Delta_X)  \qquad \text{ for } Y \in \Omega \setminus B_{\R^n}(x,2|t|).
\end{equation}

The proof of the lemma is then pretty easy. Let $Y_0$ and $X$ be as in the assumptions of the lemma. For any $Y$ far enough from $Y_0$, we use \eqref{compGw} to obtain
\[ \frac{G^Y(X)}{G^Y(Y_0)} \approx \frac{|t|^{1-d} \omega^Y(\Delta_X) }{|t_0|^{1-d} \omega^Y(\Delta_{Y_0})} ,\]
but, since the conditions on $X$ and $Y_0$ imply $E:= \Delta_X \subset 4\Delta_{Y_0}$, the estimate \eqref{chgpole} yields that
\[ \frac{G^Y(X)}{G^Y(Y_0)} \approx \left(\frac{|t|}{|t_0|}\right)^{1-d} \omega^{Y_0}(\Delta_X).\]
The above bounds on $G^Y/G^Y(Y_0)$ are uniform in $Y$, therefore, by \eqref{defG}, those bounds are transferred to $G_{Y_0}$. The lemma follows.
\ep

\section{$x$-independent Green functions with pole at infinity} \label{S3}

In this section, we shall make two easy observations: first, the Green function (associated to $L^* = -\div |t|^{n-d-1} \A^T \nabla$ as in Section \ref{SG}) with pole at infinity is independent of $x$ whenever $\A$ is $x$-independent; and second, if both $\A$ and the Green function $G$ with pole at infinity are $x$-independent, then $G$ does not depend on the first $n-d$ lines of $\A$.
We shall invoke on the top of it the uniqueness of the Green function with pole at infinity and \eqref{chgpole2} to deduce that the elliptic measure and the Lebesgue measure are equivalent on $\R^d$ whenever the last $n-d$ lines of $\A$ are $x$-independent.

\begin{lemma} \label{lemG=t}
Let $L = - \div |t|^{d+1-n} \A \nabla$ be an elliptic operator satisfying \eqref{defelliptic}--\eqref{defbounded}, and where $\A$ is as \eqref{A=AT}. Then the Green function (associated to $L^*$) with pole at infinity is $x$-independent and satisfies, for any $Y_0 = (y_0,t_0)$ and $X = (x,t) $ in $\Omega$,
\begin{equation} \label{G=t}
 C^{-1} \frac{|t|}{|t_0|} \leq G_{Y_0}(X) \leq C\frac{|t|}{|t_0|},
\end{equation} 
where the constants depend only on $n-d$ and $C_L$. 
\end{lemma}

\bp
{\bf Case 1: $(\T_3,\T_4)$ satisfies \eqref{H1}.} Define
\begin{equation} \label{defL0}
L_0 := -\div |t|^{d+1-n} \T_4 \nabla.
\end{equation}
The operator $L_0$ is an elliptic operator on $\R^{n-d} \setminus \{0\}$ satisfying the elliptic and boundedness conditions \eqref{defelliptic}--\eqref{defbounded} with the same constant $C_L$ as $L$. 

When $d<n-1$, the space $\R^{n-d} \setminus \{0\}$ and the operator $L_0$ enter the scope of the elliptic theory developed in \cite{DFMprelim} or \cite{DFMprelim2}
\footnote{Maybe also when $d=n-1$, but let us not take any risk.}, 
and so all the results in Section \ref{SG} hold. 
Of course, the study the elliptic measure of $L_0$, where the boundary is reduced to the point $\{0\}$, is trivial and hence not very interesting. But for us, this will be a chance. Let $\omega_{L_0}^X$ be the elliptic measure on $\{0\}$ and $G_{(L_0)^*,t_0}$ be the Green function with pole at infinity (associated to $(L_0)^*$) which takes the value $1$ at $t_0$. Lemma \ref{lemG=tw} implies, for $|t|<2|t_0|$, that
\[G_{(L_0)^*,t_0}(t) \approx \frac{|t|}{|t_0|} \omega^{Y_0}(\Delta_t) = \frac{|t|}{|t_0|} \omega^{Y_0}(\{0\}) = \frac{|t|}{|t_0|}. \]
The probability measure $\omega_{L_0}^X$ on $\{0\}$ obviously satisfies $\omega_{L_0}^X(\{0\}) = 1$, hence 
\begin{equation} \label{G0=t2}
G_{(L_0)^*,t_0}(t) \approx \frac{|t|}{|t_0|} \qquad \text{ for } |t|<2|t_0|. 
\end{equation}
When $|t|\geq2|t_0|$, we use \eqref{Gunique} and \eqref{G0=t2} to write 
\[G_{(L_0)^*,t_0}(t) = [G_{(L_0)^*,t}(t_0)]^{-1} \approx \left(\frac{|t_0|}{|t|} \right)^{-1} = \frac{|t|}{|t_0|}.\]
We conclude, for any $t,t_0 \in \R^{n-d}$, that
\begin{equation} \label{G0=t}
G_{(L_0)^*,t_0}(t) \approx \frac{|t|}{|t_0|}. 
\end{equation}

When $d=n-1$, the result \eqref{G0=t} holds without the need of Lemma \ref{lemG=tw}. The operator $L_0$ is defined on the half line, and there exists $f(t)$ defined on $(0,+\infty)$ such that $L_0 = \partial_t f(t) \partial_t$ and $f(t) \approx 1$ in order to satisfy the elliptic and boundedness condition. A simple exercise of integration shows that the Green functions with pole at infinity of $(L_0)^* = L_0$ are
\begin{equation} \label{zzzz}
G_{(L_0)^*}(t) = K \int_0^t \frac{dt}{f(t)} \approx C|t|,
\end{equation}
where $K$ is any positive constant, and thus \eqref{G0=t} follows easily.

\medskip

We set, for $Y_0 = (y_0,t_0) \in \Omega$ and $X= (x,t) \in \Omega$,
\begin{equation} \label{defH}
H_{Y_0}(X) := G_{(L_0)^*,t_0}(t).
\end{equation}
Check that the function $H_{Y_0}(.)$ is a Green function (associated to $L^*$) with pole at infinity and satisfies $H_{Y_0}(Y_0) = 1$. By the uniqueness given in Proposition \ref{propGinfty}, we necessary have
\begin{equation} \label{defH2}
G_{Y_0}(X) = H_{Y_0}(X) := G_{(L_0)^*,t_0}(t).
\end{equation}
As a consequence, $G_{Y_0}$ is $x$-independent, and the conclusion \eqref{G=t} of the lemma follows from \eqref{G0=t}.

\medskip 

\noindent {\bf Case 2: $(\T_3,\T_4)$ satisfies \eqref{H2}.} In this case, the proof is a simple exercise of integration. By \eqref{defelliptic} and \eqref{defbounded}, we have
\[(C_{L})^{-1} |\nabla |t||^2 \leq \T_4 \nabla |t| \cdot \nabla |t| \leq C_L |\nabla |t||^2 \qquad \text{ for all } t \in \R^{n-d}\setminus \{0\}.\]
Since $|\nabla |t|| = 1$, our assumption on $\T_4$ implies that
\[(C_{L})^{-1} \leq h(|t|) \leq C_L \qquad \text{ for all } t \in \R^{n-d}\setminus \{0\}.\]
We define $g_{r_0}$ as
\begin{equation} \label{defgr0}
g_{r_0} = K_{r_0} \int_{0}^r \frac{1}{h(r)} dr ,
\end{equation}
where $K$ is chosen so that $g_{r_0}(r_0) = 1$. Our bounds on $h$ yield that
\begin{equation} \label{gr=r}
g_{r_0} \approx \frac{r}{r_0}.
\end{equation}

We construct now $H_{Y_0}(X)$ for $Y_0 = (y_0,t_0) \in \Omega$ and $X=(x,t)\in \Omega$ as
\begin{equation} \label{defHY0}
H_{Y_0}(X) := g_{|t_0|}(|t|).
\end{equation}
Check that the definition of $g_{r_0}(r)$ and our assumption on $\A$ entail that $LH_{Y_0}(X) = 0$ and $H_{Y_0}(Y_0) = 1$, so by uniqueness of the Green function with pole at infinity (see Proposition \ref{propGinfty}), we have $G_{Y_0} = H_{Y_0}$. The conclusion \eqref{G=t} is then an easy consequence of \eqref{defHY0} and \eqref{gr=r}. \ep

\begin{remark}
An interesting consequence of the above proof, for instance \eqref{zzzz}, is that for a general operator in the form $L= -\div |t|^{d+1-n} \A \nabla$, knowing that the elliptic measure is $A_\infty$-absolute continuous with respect to the Lebesgue measure (or even equivalent to the Lebesgue measure) will not help us to get any control on $t$-derivative of the Green functions with pole at infinity. 
\end{remark}

\begin{corollary} \label{lemTt=t}
Let $L = - \div |t|^{d+1-n} \A \nabla$ be an elliptic operator satisfying \eqref{defelliptic}--\eqref{defbounded}. Assume that $\A$ can be written as
\begin{equation} 
\A = \begin{pmatrix} \A_1 & \A_2 \\ 0 & \T_4 \end{pmatrix},
\end{equation}
where $\T_4 \nabla |t| = \nabla |t|$ for all $t\in \R^{n-d}$.

Then, for $X=(x,t) \in \Omega$ and $Y_0=(y_0,t_0)\in \Omega$, the Green function with pole at infinity satisfies
\begin{equation} \label{G=t3}
G_{Y_0}(X) = \frac{|t|}{|t_0|}.
\end{equation} 
\end{corollary}

\bp Under our assumptions, $\T_4$ satisfies \eqref{H2} with $h(r) \equiv 1$. From the proof of Lemma \ref{lemG=t}, we have $G_{Y_0}(X) = g_{|t_0|}(|t|)$ where $g_{r_0}(r)$ is given by \eqref{defgr0}. The lemma follows.
\ep

We are now ready to prove Theorem \ref{Main2}.

\noindent {\bf Proof of Theorem \ref{Main2}:} Take $x\in \Delta_{Y_0}$. The combination of Lemma \ref{lemG=tw} and Lemma \ref{lemG=t} entails, for any $0<r<|t_0|$ and any $X = (x,t)$ satisfying $|t| = r$, that 
\begin{equation}
\omega^{Y_0}(B_{\R^d}(x,r)) \approx G_{Y_0}(X) \left( \frac{|t|}{|t_0|} \right)^{d-1} \approx  \left( \frac{|t|}{|t_0|} \right)^{d} = \frac{|B_{\R^d}(x,r)|}{|\Delta_{Y_0}|}.
\end{equation}
In particular, the measure is absolutely continuous with respect to the $d$-dimensional Lebesgue measure on $\R^d$ and, by the Lebesgue differentiation theorem, the Poisson kernel $k^{Y_0} := \frac{d\omega^{Y_0}}{d\mathcal L^d}$ satisfies, for almost any $x\in \Delta_{Y_0}$, that
\[ k^{Y_0}(x) = \lim_{r \to 0} \frac{\omega^{Y_0}(B_{\R^d}(x,r))}{|B_{\R^d}(x,r)|} \approx \frac{1}{|\Delta_{Y_0}|}.\]
The lemma follows by integrating $k^{Y_0}$ over $E$.
\ep

\section{Proof of Theorems \ref{Main1} and \ref{Main3}} \label{S4}

Let $1\leq d<n$ be integers, and let $\Omega = \R^{n}_+ := \{(x,t) \in \R^d \times (0,+\infty)$ if $d=n-1$ and $\Omega = \R^n \setminus \R^d:= \{(x,t) \in \R^d \times \R^{n-d}, \, t \neq 0\}$ if $d<n-1$. The non-tangential maximal functions $N$ and $\wt N$ are defined for any continuous function $v$ on $\Omega$ and any $x\in \R^d$ as
\begin{equation}
N(v)(x) = \sup_{(y,t) \in \gamma(x)} |v| \quad \text{ and } \quad \wt N(v)(x) = \sup_{(y,t) \in \gamma(x)} \left( \fint_{|Z-(y,t)| < |t|/4} |v|^2 dZ \right)^\frac12 
\end{equation}
where
\begin{equation}
\gamma(x) = \{(y,t) \in \Omega, \, |y-x|< |t|\}. 
\end{equation}
We shall introduce here the variants $\gamma_{10}(x) := \{(y,t) \in \Omega, \, |y-x| < 10|t|\}$ and $N_{10}(v)(x) := \sup_{\gamma_{10}(x)} |v|$. They will be used to compare $\wt N$ and $N$. Indeed, we have the pointwise bound $\wt N(v)(x) \leq N_{10}(v)(x)$ and it is well known (see \cite{Stein93}, Chapter II, \S 2.5.1) that $\|N_{10}(v)\|_{2} \approx \|N(v)\|_{2}$. Altogether, 
\begin{equation} \label{wtNN}
\|\wt N(v)\|_{L^2(\R^d)} \leq \|N_{10}(v)\|_{L^2(\R^d)} \approx \|N(v)\|_{L^2(\R^d)}.
\end{equation}

We recall that the non-tangential maximal functions behave well with the Carleson measure condition \eqref{defCM} and \eqref{defwtCM}. Indeed, if $v$ is a continuous function on $\Omega$ and $f\in CM_2(K)$, then we have the Carleson inequality
\begin{equation} \label{Carleson}
\int_\Omega f^2 v^2 \frac{dx \, dt}{|t|^{n-d}} \lesssim K \|N(v)\|_{L^2(\R^d)}^2, 
\end{equation}
and similarly, if $g\in \wt{CM}_2(K)$, then 
\begin{equation} \label{Carleson2}
\int_\Omega g^2 v^2 \frac{dx \, dt}{|t|^{n-d}} \lesssim K \|\wt N(v)\|_{L^2(\R^d)}^2 \lesssim K \|N(v)\|_{L^2(\R^d)}^2.
\end{equation}
Combined with the Cauchy-Schwarz inequality, for all $w\in L^2_{loc}(\Omega)$, one has
\begin{equation} \label{CS+C}
\int_\Omega f v w \frac{dx \, dt}{|t|^{n-d}} \leq C K^{\frac12} \|N(v)\|_{L^2(\R^d)} \left(\int_\Omega w^2 \frac{dx \, dt}{|t|^{n-d}}\right)^\frac12 
\end{equation}
and 
\begin{equation} \label{CS+C2}
\int_\Omega g v w \frac{dx \, dt}{|t|^{n-d}} \leq C K^{\frac12} \|\wt N(v)\|_{L^2(\R^d)} \left(\int_\Omega w^2 \frac{dx \, dt}{|t|^{n-d}}\right)^\frac12.
\end{equation}

\medskip

We also introduce cut-off functions associated to tent sets. Choose a smooth function $\phi \in C^\infty_0(\R)$ such that $0 \leq \phi \leq 1$, $\phi\equiv 1$ on $(-1,1)$, $\phi \equiv 0$ on $(2,+\infty)$ and $|\phi'|\leq 2$. For a ball $B := B(x,r) \subset \R^d$, we define $\Psi_B$ as
\begin{equation} \label{defPsiB}
\Psi_B(y,t) = \phi\left( \frac{\dist(x,B)}{|t|}\right) \phi\left(\frac{|t|}{r}\right).
\end{equation}
We also associate to $B$ the tent set $T_B:= \{(x,t) \in \Omega, \, x\in B, \, |t|\leq r \}$. The function $\Psi_B$ is such that $\Psi \equiv 1$ on $T_B$ and $\Psi \equiv 0$ on $\Omega \subset T_{2B}$. Note that, if a different definition of tent sets is used, we can easily change the definition of $\Psi_B$ so that $\Psi_B$ is adapted to the other definition of tent sets.

Theorems \ref{Main1} and \ref{Main3} are consequences of the following lemma.

\begin{lemma}
If $L = - \div |t|^{d+1-n} \nabla$ satisfies the assumptions of Theorem \ref{Main1} or \ref{Main3}, then, for any ball $B=B(x,r) \subset \R^d$, we have
\begin{equation} \label{S<N}
\int_\Omega |\nabla u|^2 \Psi_B^4 \frac{dt \, dx}{|t|^{n-d-2}} \leq C (1+ K) \|N(u\Psi_B)\|_{L^2(\R^d)}^2,
\end{equation}
where $C>0$ depends only on $n$, $d$, and $C_L$.
\end{lemma}

\bp {\bf Step 1: Carleson estimates on the cut-off functions.}  In order to deal with finite quantities, we need to refine our cut-off function $\Psi_{B}$. We define $\Psi_{B,\epsilon}$ as 
\begin{equation} \label{defPsiBe}
\Psi_{B,\epsilon} (y,t) = \Psi_B(y,t) \phi \left( \frac{\epsilon}{|t|} \right),
\end{equation}
where $\phi$ is the smooth function introduced above \eqref{defPsiB} and already used to define $\Psi_B$ . We first gather some properties of the cut-off function $\Psi_{B,\epsilon}$. Observe that
\begin{equation} \label{defPsiBeZ}
|\nabla \Psi_{B,\epsilon}(y,t)| \lesssim \frac{1}{|t|} \quad \text{ for } (y,t) \in \Omega,
\end{equation}
and $\nabla \Psi_{B,\epsilon}$ is supported on $E_1 \cup E_2 \cup E_3$, where
\[E_1:= \{(y,t)\in \Omega, \, \dist(y,B) \leq 2|t| \leq 2\dist(y,B) \},\]
\[E_2:= \{(y,t)\in \Omega, \, r(B) \leq |t| \leq 2r(B)\},\]
with $r(B)$ being the radius of $B$, and
\[E_3:= \{(y,t)\in \Omega, \, |t| \leq \epsilon \leq 2|t|\}.\]
 So we deduce that 
 \begin{equation} \label{defPsiBeY}
|t||\nabla \Psi_{B,\epsilon}(y,t)| +  |t|^2|\nabla \Psi_{B,\epsilon}(y,t)|^2 \lesssim \1_{E_1 \cup E_2 \cup E_3}(y,t).
\end{equation}
We will need the fact that $|t||\nabla \Psi_{B,\epsilon}(y,t)|$ and $(|t||\nabla \Psi_{B,\epsilon}(y,t)|)^{1/2}$ satisfy the Carleson measure condition $\wt{CM}_2(M)$ for some uniform constant $M$ which, combined with \eqref{Carleson}, implies, for any continuous functions $v$, that
 \begin{equation} \label{tNablaPsiCM}
\int_\Omega |t| |\nabla \Psi_{B,\epsilon}(y,t)| v^2 \frac{dt\, dx}{|t|^{n-d}} + \int_\Omega |t| |\nabla \Psi_{B,\epsilon}(y,t)| v^2 \frac{dt\, dx}{|t|^{n-d}} \lesssim \|\wt N(v)\|_{L^2(\R^d)}^2 .
\end{equation}
Of course, thanks to \eqref{wtNN}, if \eqref{tNablaPsiCM} is true, then we also have the analogue estimate where $\wt N$ is replaced by $N$. Thanks to \eqref{defPsiBeY}, the claim \eqref{tNablaPsiCM} will be then proven if we can show that $\1_{E_1 \cup E_2 \cup E_3} \in \wt{CM}_2(M)$, that is
 \begin{equation} \label{EiCM}
\sup_{x\in \R^d, \, r>0} \fint_{B_{\R^d}(x,r)} \int_{|t| < r} \sup_{|Z-(y,t)|< |t|/4} |\1_{E_1 \cup E_2 \cup E_3}(Z)|^2 \frac{dy \, dt}{|t|^{n-d}} \lesssim 1.
\end{equation}
However, \eqref{EiCM} is an immediate consequence of the fact that, for each $y\in \R^d$, there holds
\[\begin{split}
\int_{t \in \R^{n-d}}  \sup_{|Z-(y,t)|< |t|/4} |\1_{E_1 \cup E_2 \cup E_3}(Z)|^2  \frac{dt}{|t|^{n-d}} 
& \leq  \int_{\dist(y,B)/4 \leq |t| \leq 2\dist(y,B)} \frac{dt}{|t|^{n-d}} \\
& \hspace{0.3cm} + \int_{r(B)/2 \leq |t| \leq 4r(B)} \frac{dt}{|t|^{n-d}} + \int_{\epsilon/4 \leq |t| \leq 2\epsilon} \frac{dt}{|t|^{n-d}} \lesssim 1.
\end{split}\]
The claim \eqref{tNablaPsiCM} follows.

\medskip

\noindent  {\bf Step 2: Introduction of $G$.}  First, we decompose $L$ as
\begin{equation} \label{A=B+C5}
\A = \begin{pmatrix} \A_1 & \A_2 \\ \B_3 + \wt \B_3 + \C_3 & b \T_4 + \C_4 \end{pmatrix} 
\end{equation}
so that it includes both the case of Theorem \ref{Main1} and \ref{Main3}. In particular, we have 
\begin{equation} \label{CMK} \begin{split}
|\C_3| + |\C_4| & \in \wt{CM}_2(K) \\
 |t| |\nabla b| + |t||\div_x (\B_3)^T|  & +  |t|^{n-d} |\diver_t (|t|^{d+1-n} \B_3)| \\
& \hspace{2cm} + |t||\div_x (\wt \B_3)^T| + |t| |\div_t \wt\B_3|   \in CM_2(3K).
\end{split} \end{equation}
We set $L_0:= - \div |t|^{d+1-n} \A_0 \nabla$ where
\begin{equation} \label{defA0}
\A_0 := \begin{pmatrix} \A_1/b & \A_2/b \\ 0 & \T_4 \end{pmatrix},
\end{equation}  
which satisfies
\begin{equation} \label{A-A0}
\frac{1}{b} \A = \A_0 + \frac1b  \begin{pmatrix} 0 & 0 \\ \B_3 + \tilde \B_3 + \C_3 & \C_4 \end{pmatrix}.
\end{equation}  
Let $Y_0 = (y_0,t_0) \in \Omega$ be such that $|t_0|= 1$, and we write $G$ for $G_{Y_0}$, the Green function associated to $(L_0)^*$ with pole at infinity.
The important properties on $G$ for this proof are first that $G \in W^{1,2}_{loc}(\Omega)$ is a weak solution to $(L_0)^* u = 0$ in $\Omega$, that is
\begin{equation} \label{Gweak}
\int_\Omega  \A_0 \nabla \varphi \cdot \nabla G \, \frac{dt \, dx}{|t|^{n-d-1}} \qquad \text{ for any compactly supported $\varphi \in W^{1,2}(\Omega)$},
\end{equation}
and second that, Lemma \ref{lemG=t} entails that
\begin{equation} \label{G=tp}
\text{$G$ is $x$-independent and }  G(X) \approx |t| \text{ for all } X = (x,t) \in \Omega.
\end{equation}

\medskip

\noindent {\bf Step 3: Estimation of $\|\wt N(u\Psi_{B,\epsilon}^2\nabla G)\|_2$.} 
If the goal were to only obtain \eqref{NablauCM}, we would not need to go through the same computations, we just have to prove
\begin{equation} \label{claimS3}
\|\wt N(u_H\Psi_{B,\epsilon}^2\nabla G)\|_{L^2(\R^d)} \lesssim |B|.
\end{equation}
Since $G$ is a weak solution to $L_0 u = 0$, Caccioppoli's inequality yields that
\[\fint_{|Z-(y,t)| < |t|/4} |\nabla G|^2 dZ \lesssim \frac1{|t|^2} \fint_{|Z-(y,t)| < |t|/2} |G|^2 dZ \qquad \text{ for } (y,t) \in \Omega.\]
 But since $G\approx |t|$ by \eqref{G=tp}, the above inequality becomes
\[ \fint_{|Z-(y,t)| < |t|/4} |\nabla G|^2 dZ \lesssim 1.\]
We take the supremum on $(y,t) \in \gamma(x)$, and then integrate on $x\in 100B$, and we get 
\[|B| \gtrsim \|\wt N(\nabla G)\|_{L^2(100B)} \gtrsim \|\wt N(u_H\Psi_{B,\epsilon}^2\nabla G)\|_{L^2(\R^d)} \]
because $u_H \leq 1$ by construction. The claim \eqref{claimS3} follows.
 
\medskip

However, what we really need in order to prove the inequality \eqref{S<N} is
\begin{equation} \label{claimS4}
\|\wt N(u\Psi_{B,\epsilon}^2\nabla G)\|_{L^2(\R^d)} \lesssim \|N(u\Psi_{B,\epsilon})\|_{L^2(\R^d)}.
\end{equation}
To reach this goal, we first need the following Caccioppoli's inequality. Let $D \subset \R^n$ be a ball of radius $r$ such that $4D \subset \Omega$ and $5D\cap \partial \Omega \neq \emptyset$. In particular, we have
\begin{equation} \label{truc2}
G(X) \approx |t| \approx r \quad \text{ for } X=(x,t) \in 2D
\end{equation}
by \eqref{G=tp}. Let $\Psi$ be a function such that $0\leq \Psi\leq 1$ and $|\nabla \Psi| \lesssim \frac1{|t|}$, and let $u$ bs a weak solution to $Lu = 0$. We claim that
\begin{equation} \label{claimS5}
\fint_{D} |\nabla G|^2 u^2 \Psi^4 \, dX \lesssim \frac1{r^2} \fint_{2D} |u|^2 \Psi^2 \, dX. 
\end{equation}
Let $\Phi$ be such that $0 \leq \Phi \leq 1$, $\Phi \equiv 1$ on $D$, $\Phi \equiv 0$ outside $4D/3$, and $|\nabla \Phi| \leq 5r$. Then
\begin{equation} \label{truc1}
\int_{D} |\nabla G|^2 u^2 \Psi^4 \, dX \leq T:= \int_{D} |\nabla G|^2 u^2 \Psi^4 \Phi^2 \, dX.
\end{equation}
The function $G$ is a weak solution of $L_0 u = 0$ so, by ellipticity of $\A_0$ and since the weight satisfies $|t|^{d+1-n} \approx r^{d+1-n}$ on $2D$, we have
\[ \begin{split}
T &  \lesssim \iint_{\Omega} \A_0 \nabla G \cdot \nabla G u^2 \Psi^4 \Phi^2 \, \frac{dt \, dx}{|t|^{n-d-1}} \\
& =  \iint_{\Omega} \A_0 \nabla [G u^2 \Psi^4 \Phi^2]  \cdot \nabla G  \, \frac{dt \, dx}{|t|^{n-d-1}} - 2 \iint_{\Omega} \A_0 \nabla u \cdot \nabla G  \, (Gu \Psi^4 \Phi^2) \, \frac{dt \, dx}{|t|^{n-d-1}} \\
& \quad - 2 \iint_{\Omega} \A_0 \nabla \Phi \cdot \nabla G \, (G u^2 \Psi^4 \Phi) \, \frac{dt \, dx}{|t|^{n-d-1}} - 4 \iint_{\Omega} \A_0 \nabla \Psi  \cdot \nabla G \, (Gu^2 \Psi^3 \Phi^2) \, \frac{dt \, dx}{|t|^{n-d-1}}\\
& := T_1 + T_2 + T_3 + T_4.
\end{split} \]
The functions $G$, $u$, $\Phi$, and $\Psi$ all belong to $L^\infty(2D) \cap W^{1,2}(2D)$, so $G u^2 \Psi^4 \Phi^2$ is a valid test function and \eqref{Gweak} gives that $T_1 = 0$. By the boundedness of $\A_0$ and Cauchy-Schwarz's inequality, the terms $T_2$, $T_3$, and $T_4$ can be bounded as follows. We have
\[ \begin{split}
|T_3| &  \lesssim T^{1/2} \left( \iint_\Omega |\nabla \Phi|^2 G^2  u^2 \Psi^4 \frac{dt \, dx}{|t|^{n-d-1}} \right)^{\frac12} \lesssim T^{1/2} \left( \iint_{4D/3} u^2 \Psi^4 \frac{dt \, dx}{|t|^{n-d-1}} \right)^{\frac12},
\end{split} \]
because $|\nabla \Phi| \lesssim 1/r \approx 1/G$ on $2D$. Similarly
\[ \begin{split}
|T_4| &  \lesssim T^{1/2} \left( \iint_\Omega  |\nabla \Psi| G^2 u^2 \Psi^2 \Phi^2  \frac{dt \, dx}{|t|^{n-d-1}} \right)^{\frac12} \lesssim T^{1/2} \left( \iint_{4D/3}  u^2 \Psi^2  \frac{dt \, dx}{|t|^{n-d-1}} \right)^{\frac12}
\end{split} \]
because $|\nabla \Psi| \lesssim \frac1{|t|} \approx \frac1{G}$ on $2D$. At last, 
\[ \begin{split}
|T_2| &  \lesssim T^{1/2} \left( \iint_{\Omega} |\nabla u|^2 G^2 \Psi^4 \Phi^2 \frac{dt \, dx}{|t|^{n-d-1}} \right)^{\frac12} \lesssim T^{1/2} \left( r^2 \iint_{4D/3} |\nabla u|^2 \Psi^4 \frac{dt \, dx}{|t|^{n-d-1}} \right)^{\frac12}.
\end{split} \]
We deduce that 
\[T \lesssim T^{1/2} \left( \iint_{4D/3} u^2 \Psi^2 \frac{dt \, dx}{|t|^{n-d-1}} + r^2 \iint_{4D/3} |\nabla u|^2 \Psi^4 \frac{dt \, dx}{|t|^{n-d-1}} \right)^\frac12\]
and then
\begin{equation} \label{truc3}
\iint_D |\nabla G|^2 u^2 \Psi^4 \, \frac{dt \, dx}{|t|^{n-d-1}} \lesssim  \iint_{4D/3} u^2 \Psi^2 \frac{dt \, dx}{|t|^{n-d-1}} + r^2 \iint_{4D/3} |\nabla u|^2 \Psi^4 \frac{dt \, dx}{|t|^{n-d-1}}.
\end{equation}
We repeat the process for the last integral of the right-hand side above, using the fact that $u$ is a weak solution to $Lu=0$, and we obtain\footnote{The estimate below can also be seen as a variant of Caccioppoli's inequality, and is a consequence of Lemma 3.1 (i) in \cite{FMZ}.}
\[r^2 \iint_{4D/3} |\nabla u|^2 \Psi^4 \frac{dt \, dx}{|t|^{n-d-1}} \lesssim \iint_{2D} u^2 \Psi^2 \frac{dt \, dx}{|t|^{n-d-1}}.\]
We combine the last estimate with \eqref{truc3}, and we get that
\begin{equation} \label{truc4}
\iint_D |\nabla G|^2 u^2 \Psi^4 \, \frac{dt \, dx}{|t|^{n-d-1}} \lesssim  \iint_{2D} u^2 \Psi^2 \frac{dt \, dx}{|t|^{n-d-1}}.
\end{equation}
The claim \eqref{claimS5} follows after we recall that $|t| \approx r$ on $2D$.

\medskip

We apply now \eqref{claimS5} to have
\[\fint_{|Z-(y,t)| < |t|/4} |\nabla G|^2 u^2 \Psi_{B,\epsilon}^4 dZ \lesssim \frac1{|t|^2} \fint_{|Z-(y,t)| < |t|/2} u^2 \Psi_{B,\epsilon}^2 dZ \qquad \text{ for } (y,t) \in \Omega.\]
As a consequence, for any $x\in \R^d$,
\[\wt N(u\Psi_{B,\epsilon}^2\nabla G)(x) \lesssim N_{10}(u\Psi_{B,\epsilon})(x).\]
The claim \eqref{claimS4} follows from \eqref{wtNN}.

\medskip

\noindent {\bf Step 4: Proof of \eqref{S<N}.}  We define 
\[J = J_{B,\epsilon}:= \int_\Omega |\nabla u|^2 \Psi_{B,\epsilon}^4 \frac{dt\, dx}{|t|^{n-d-2}}\]
and we want to show that 
\begin{equation} \label{defJBe}
J_{B,\epsilon} \lesssim (1+K) \|N(u\Psi_{B,\epsilon} )\|_{L^2(\R^d)}^2 + (1+K^{1/2}) J_{B,\epsilon}^{1/2} \|N(u\Psi_{B,\epsilon} )\|_{L^2(\R^d)},
\end{equation}
where $K$ is the constant used in the assumptions of the theorem under proof.
Since $u\in W^{1,2}_{loc}(\Omega)$, all the quantities in \eqref{defJBe} are finite, and therefore \eqref{defJBe} improves itself in
\begin{equation} \label{defJBe2}
J_{B,\epsilon} \lesssim (1+K) \|N(u\Psi_{B,\epsilon} )\|_{L^2(\R^d)}^2.
\end{equation}
We assumed that the solution $u$ is bounded, so the left-hand side above is uniformly bounded in $\epsilon$. We take then the limit as $\epsilon$ goes to 0 to obtain the desired bound \eqref{S<N}.

\medskip

To lighten the notation, we shall write until the end of the proof $\Psi$ for $\Psi_{B,\epsilon}$ and $J$ for $J_{B,\epsilon}$. 
Since $b$ is bounded from above (assumption $(b)$ of both Theorems \ref{Main1} and \ref{Main3}), $G \gtrsim |t|$ by \eqref{G=tp}, and $\A$ is elliptic by \eqref{defelliptic}, we deduce that
\[J \lesssim I:= \iint_\Omega \A \nabla u \cdot \nabla u \frac{\Psi^4 G}{b} \frac{dt\, dy}{|t|^{n-d-1}}.\]
Using the product rule, we insert $\Psi^4 G/b$ into the second gradient, and we obtain
\[\begin{split}
I & = \iint_\Omega \A \nabla u \cdot \nabla \left( \frac{u \Psi^4 G}{b}\right) \frac{dt\, dy}{|t|^{n-d-1}} - 4 \iint_\Omega \A \nabla u \cdot \nabla \Psi \, \frac{u \Psi^3 G}{b} \frac{dt\, dy}{|t|^{n-d-1}} \\
& \quad + \iint_\Omega \A \nabla u \cdot \nabla b \, \frac{u\Psi^4 G}{b^2} \frac{dt\, dy}{|t|^{n-d-1}} - \iint_\Omega \A \nabla u \cdot \nabla G \, \frac{u\Psi^4}{b} \frac{dt\, dy}{|t|^{n-d-1}} \\
& = I_0 + I_1 + I_2 + I_3
\end{split}\]
The term $I_0$ equals $0$ because $u$ is a weak solution to $Lu=0$ (and the compactly supported function $u \Psi^4 G/b \in W^{1,2}(\Omega)$ is a valid test function thanks to Lemma 8.3 in \cite{DFMprelim}). 
The terms $I_1$ and $I_2$ are bounded in a similar manner. Since $b \gtrsim 1$, $G\approx |t|$, $\A$ is bounded (due to \eqref{defbounded}), and $0\leq \Psi \leq 1$, the Cauchy-Schwarz inequality infers that
\[\begin{split}
|I_1+I_2| & \lesssim \iint_\Omega |t|(|\nabla \Psi| + |\nabla b|) u \Psi^3 |\nabla u| \frac{dt\, dy}{|t|^{n-d-1}} \\
& \lesssim J^{1/2} \left(  \iint_\Omega |t|^2(|\nabla \Psi|^2 + |\nabla b|^2) u^2 \Psi^2 \frac{dt\, dy}{|t|^{n-d}} \right)^\frac12.
\end{split}\]
We know that $|t||\nabla b| \in CM(K)$ by assumption \eqref{CMK}, and that $|t||\nabla \Psi| \in CM$ by \eqref{tNablaPsiCM}, so the Carleson inequality \eqref{Carleson} entails that
\[|I_2 + I_3| \lesssim (1+K^{1/2}) J^{1/2} \|N(u\Psi)\|_{L^2(\R^d)}.\] 

As for $I_3$, we use the decomposition of $\A$ given in \eqref{A-A0} to obtain 
\[\begin{split}
I_3 & = - \iint_\Omega (\C_3 \nabla_x u + \C_4\nabla_t u) \cdot \nabla_t G \, \frac{u\Psi^4}{b} \frac{dt\, dy}{|t|^{n-d-1}} - \iint_\Omega  \A_0 \nabla u \cdot \nabla G \, (u\Psi^4) \frac{dt\, dy}{|t|^{n-d-1}}   \\
& \hspace{2cm}  - \iint_\Omega (\B_3+\wt \B_3) \nabla_x u \cdot \nabla_t G \, \frac{u\Psi^4}{b} \frac{dt\, dy}{|t|^{n-d-1}} \\
& = I_{31} + I_{32} + I_{33}.
\end{split}\] 
Recall that $b\gtrsim 1$, combined with the fact that  $|\C_3|+ |\C_4| \in \wt{CM}_2(K)$ and \eqref{CS+C2}, we deduce
\[\begin{split}
|I_{31}| & \lesssim \iint_\Omega (|\C_3|+ |\C_4|) |\nabla u| |u| |\nabla G| \Psi^4 \frac{dt\, dy}{|t|^{n-d-1}} \lesssim J^{1/2} K^{1/2} \|\wt N(u\Psi^2 \nabla G)\|_{L^2(\R^d)} \\
& \lesssim J^{1/2} K^{1/2} \|N(u\Psi)\|_{L^2(\R^d)}
\end{split}\]
by \eqref{CS+C} and then \eqref{claimS4}. We force $(u\Psi^4)$ into the first gradient and $T_{32}$ becomes
\[\begin{split}
I_{32} &  = - \frac12 \iint_\Omega \A_0 \nabla (u^2 \Psi^4) \cdot \nabla G \frac{dt\, dy}{|t|^{n-d-1}}
) + 2 \iint_\Omega \A_0 \nabla \Psi \cdot \nabla G \, (u^2\Psi^3) \frac{dt\, dy}{|t|^{n-d-1}} \\
& := I_{321} + I_{322}  \end{split}\]
The term $I_{321} = 0$ thanks to \eqref{Gweak}. As for $I_{322}$, we use the boundedness of $\A_0$ and the inequality that $2ab \leq a^2+b^2$ to write
\[I_{322} \lesssim \iint_\Omega |\nabla \Psi|  |\nabla G|^2 u^2\Psi^4 \, \frac{dt\, dy}{|t|^{n-d-1}} + \iint_\Omega |\nabla \Psi|  u^2\Psi^2 \frac{dt\, dy}{|t|^{n-d-1}},\]
and then, by \eqref{tNablaPsiCM} and \eqref{claimS4},
\[I_{322} \lesssim \|\wt N(u\Psi^2\nabla G)\|_{L^2(\R^d)}^2 + \|\wt N(u\Psi)\|_{L^2(\R^d)}^2 \lesssim \|N(u\Psi)\|_{L^2(\R^d)}^2.\]

\medskip

\noindent {\bf Step 5: Bound of $I_{33}$, which is the only difference between Theorem \ref{Main1} and \ref{Main3}.} The goal is to permute the two gradients $\nabla_x$ and $\nabla_t$ on $I_{33}$. We define $S$ as the part of $I_{33}$ that contains $\B_3$, i.e.
\begin{equation} \label{defSj}
S := - \iint_\Omega \B_3 \nabla_x u \cdot \nabla_t G \, \frac{u\Psi^4}{b} \, \frac{dt\, dy}{|t|^{n-d-1}}.
\end{equation}
Using integration by parts in $t$, $S$ becomes
\[\begin{split}
S & = -\frac12 \int_\Omega \B_3\nabla_x [u^2] \cdot \nabla_t G \,  \frac{\Psi^4}{b} \,  \frac{dt\, dy}{|t|^{n-d-1}} \\
& = \frac12 \int_\Omega \div_t(|t|^{d+1-n} \B_3\nabla_x [u^2])  \frac{G\Psi^4}{b} \, dt\, dy + 2 \int_\Omega \B_3 \nabla_x[u^2] \cdot \nabla_t \Psi \, \frac{G\Psi^3}{b}  \frac{dt\, dy}{|t|^{n-d-1}} \\
& \qquad - \int_\Omega \B_3 \nabla_x[u^2] \cdot \nabla_t b \, \frac{G\Psi^4}{b^2}  \frac{dt\, dy}{|t|^{n-d-1}}  \\
& : = S_0 + S_1 + S_2
\end{split}\]
We write the term $S_0$ as a sum on the coefficients of $\B_3$, we permute the $x$ and the $t$-derivatives on $u^2$, and then we integrate by parts in $x$. Recall that, in this paper, when $M$ is a matrix-valued functions, $\div M$ are vector-valued functions whose $j^{th}$ entry is the divergence of the $j^{th}$ column of $M$. 
\[\begin{split}
S_0 & := \frac12 \sum_{1\leq j \leq d < i \leq n}  \iint_\Omega \partial_{t_i} [ |t|^{d+1-n} (\B_3)_{ij} \partial_{x_j} u^2]  \frac{G\Psi^4}{b} \, dt\, dy \\
&  = \frac12 \iint_\Omega \diver_t(|t|^{d+1-n} \B_3) \cdot \nabla_x (u^2) \frac{G\Psi^4}{b}\, dt\, dy + \frac12 \sum_{1\leq j \leq d < i \leq n}  \iint_\Omega (\B_3)_{ij}  \partial_{t_i} [\partial_{x_j} u^2]  \frac{G\Psi^4}{b} \frac{dt\, dy}{|t|^{n-d-1}} \\
&= \frac12 \iint_\Omega \diver_t(|t|^{d+1-n} \B_3) \cdot \nabla_x (u^2) \frac{G\Psi^4}{b}  \, dt\, dy - \frac{1}2  \iint_\Omega \diver_x (\B_3)^T \cdot \nabla_t [ u^2]  \frac{G\Psi^4}{b} \frac{dt\, dy}{|t|^{n-d-1}} \\
& \qquad -2  \iint_\Omega (\B_3)^T \nabla_t [ u^2] \cdot \nabla_x \Psi \, \frac{G\Psi^3}{b} \frac{dt\, dy}{|t|^{n-d-1}} + \frac12 \iint_\Omega (\B_3)^T \nabla_t [ u^2] \cdot \nabla_x b \,  \frac{G\Psi^3}{b^2} \frac{dt\, dy}{|t|^{n-d-1}} \\
& := S_3 + S_4 + S_5 + S_6.
\end{split}\]
We do not have $x$-derivatives on $G$ because $G$ is $x$-independent, see \eqref{G=tp}.
We deal with $S_1$, $S_2$, $S_3$, $S_4$, $S_5$, and $S_6$ in a similar manner as $I_2+I_3$ earlier. We have $G \lesssim |t|$, $1/b \lesssim 1$, and $\B_3$ is bounded, hence, if 
\[f := |t||\nabla \Psi| + |t||\nabla b| + |t||\diver_x (\B_3)^T| + |t|^{n-d} |\diver_t(|t|^{d+1-n} \B_3),\]
the sum of the $S_i$ can be bounded by
\[\begin{split} 
|S| \leq \sum_{i=1}^6 |S_i| & \lesssim \iint_\Omega f |\nabla (u^2)| \Psi^3 \frac{dt\, dy}{|t|^{n-d-1}} \\
& \lesssim \iint_\Omega f |\nabla u| u \Psi^3 \frac{dt\, dy}{|t|^{n-d-1}}.
\end{split}\]
But, since $f\in CM_2(1+K)$ by \eqref{CMK} and \eqref{tNablaPsiCM}, the Carleson estimate \eqref{CS+C} yields that 
\begin{equation} \label{boundS}
|S| \leq \sum_{i=1}^6 |S_i| \lesssim J^{1/2} (1+K^{1/2}) \|N(u\Psi)\|_{L^2(\R^d)}
\end{equation}
as desired. Theorem \ref{Main1} is now proven, because $\wt \B_3 = 0$ in its assumption.

\medskip

In order to establish Theorem \ref{Main3}, it remains to treat the part of $I_{33}$ that contains $\wt \B_3$. If $\wt S := I_{33} - S$, and if, for $i\in \{0,\dots,6\}$, $\wt S_i$ is obtained from $S_i$ by substituting $\B_3$ for $\wt \B_3$. For $i\neq 3$, we can bound $\wt S_i$ as we bound $S_i$, because the assumptions on $\wt \B_3$ match those of $\B_3$. So, similarly to \eqref{boundS}, we have that
\begin{equation} \label{boundwtS}
|\wt S - \wt S_3| \lesssim J^{1/2} (1+K^{1/2}) \|N(u\Psi)\|_{L^2(\R^d)}.
\end{equation}
We do not know how to estimate $\wt S_3$, but instead we know how to estimate
\begin{equation} \label{defS7}
\wt S_7 :=  \frac12 \iint_\Omega \diver_t(\wt \B_3) \cdot \nabla_x (u^2) \frac{G\Psi^4}{b} \frac{dt\, dy}{|t|^{n-d-1}}  = \iint_\Omega \diver_t(\wt \B_3) \cdot \nabla_x u \frac{Gu\Psi^4}{b} \frac{dt\, dy}{|t|^{n-d-1}}.
\end{equation}
Indeed, we use $G\lesssim |t|$, $1/b \lesssim 1$, $|t||\diver_t(\B_3)| \in CM_2(K)$, and the Carleson estimate \eqref{CS+C}, to get, similarly to the $S_i$'s, that 
\begin{equation} \label{boundS7}
|\wt S_7| \lesssim J^{1/2} (1+K^{1/2}) \|N(u\Psi)\|_{L^2(\R^d)}.
\end{equation}
So, in order to bound $\wt S$ and prove Theorem \ref{Main3}, we only have to write $\wt S$ as a linear combination of $\wt S_7$ and $\wt S_3$. Since we are currently under the assumptions of Theorem \ref{Main3}, Lemma \ref{lemTt=t} entails that $G = |t|$. With this in mind, we have
\[ G|t|^{n-d-1} \diver_t(|t|^{d+1-n}\B_3) = G \diver_t (\B_3) + (d+1-n) (\nabla_t G)^T \B_3,\]
which can be reformulated as
\[\wt S_3 = \wt S_7 + (n-d-1) \wt S.\]
We conclude that 
\[|\wt S| = \frac1{n-d-2} |(\wt S_3 - \wt S) + \wt S_7 | \lesssim J^{1/2} (1+K^{1/2}) \|N(u\Psi)\|_{L^2(\R^d)}\]
by \eqref{boundwtS} and \eqref{boundS7}. The lemma follows.
\ep

\end{document}